# AN EFFICIENT AND FAST SPARSE GRID ALGORITHM FOR HIGH-DIMENSIONAL NUMERICAL INTEGRATION

HUICONG ZHONG[*] AND XIAOBING FENG[†]

**Abstract.** This paper is concerned with developing an efficient numerical algorithm for fast implementation of the sparse grid method for computing the $d$-dimensional integral of a given function. The new algorithm, called the MDI-SG (*multilevel dimension iteration sparse grid*) method, implements the sparse grid method based on a dimension iteration/reduction procedure, it does not need to store the integration points, neither does it compute the function values independently at each integration point, instead, it re-uses the computation for function evaluations as much as possible by performing the function evaluations at all integration points in cluster and iteratively along coordinate directions. It is showed numerically that the computational complexity (in terms of CPU time) of the proposed MDI-SG method is of polynomial order $O(Nd^3)$ or better, compared to the exponential order $O(N(\log N)^{d-1})$ for the standard sparse grid method, where $N$ denotes the maximum number of integration points in each coordinate direction. As a result, the proposed MDI-SG method effectively circumvents the curse of dimensionality suffered by the standard sparse grid method for high-dimensional numerical integration.

**1. Introduction.** With rapid developments in nontraditional applied sciences such as mathematical finance [11], image processing [15], economics [1], and data science [24], there is an ever increasing demand for efficient numerical methods for computing high-dimensional integration which also becomes crucial for solving some challenging problems. Numerical methods (or quadrature rules) mostly stem from approximating the Riemann sum in the definition of integrals, hence, they are grid-based. The simplest and most natural approach for constructing numerical quadrature rules in high dimensions is to apply the same 1-d rule in each coordinate direction, this then leads to tensor-product (TP) quadrature rules. It is well known (and easy to check) that the number of integration points (and function evaluations) grows *exponentially* in the dimension $d$, such a phenomenon is known as *the curse of dimensionality* (CoD). Mitigating or circumventing the CoD has been the primary goal when it comes to constructing efficient high-dimensional numerical quadrature rules. A lot of progress has been made in this direction in the past fifty years, this includes sparse grid (SG) methods [4, 10, 11, 7], Monte Carlo (MC) methods [5, 21], Quasi-Monte Carlo (QMC) methods [6, 13, 14, 16, 27], deep neural network (DNN) methods [8, 12, 17, 25, 28]. To some certain extent, those methods are effective for computing integrals in low and medium dimensions (i.e., $d \lesssim 100$), but it is still a challenge for them to compute integrals in very high dimensions (i.e., $d \approx 1000$).

This is the second installment in a sequel [9] which aims at developing fast numerical algorithms for high-dimensional numerical integration. As mentioned above, the straightforward implementation of the TP method will evidently run into the CoD dilemma. To circumvent the difficulty, we proposed in [9] a multilevel dimension iteration algorithm (called MDI-TP) for accelerating the TP method. The ideas of the MDI-TP algorithm are to reuse the computation of function evaluation as much as possible in the tensor product method by clustering computations, which allows an efficient and fast function evaluations at integration points together, and to do clustering by a simple dimension iteration/reduction strategy, which is possible because of the lattice structure of the TP integration points. Since the idea of the MDI strategy essentially applied to any numerical integration rule whose integration points have a lattice-like structure, this indeed motivates the work of this paper by applying the MDI idea to accelerate the sparse grid method.

---

*School of Mathematics and Statistics, Northwestern Polytechnical University, Xi'an, Shaanxi 710129, China (huicongzhong@mail.nwpu.edu.cn).

†Department of Mathematics, The University of Tennessee, Knoxville, TN, 37996 (xfeng@utk.edu).

The sparse grid (SG) method, which was first proposed by Smolyak in [26], only uses a (small) subset of the TP integration points while still maintains a comparable accuracy of the TP method. As mentioned earlier, the SG method was one of few successful numerical methods which can mitigate the CoD in high-dimensional computation, including computing high-dimensional integration and solving high-dimensional PDEs [4, 10, 11, 7]. The basic idea in application to high-dimensional numerical integration stems from Smolyak's general method for multivariate extensions of univariate operators. Based on this construction, the midpoint rule [2], the rectangle rule [22], the trapezoidal rule [3], the Clenshaw-Curtis rule [19, 20], and the Gaussian-Legendre rule [18, 23] have been used as a one-dimensional numerical integration method. A multivariate quadrature rule is then constructed by forming the TP method of each of these one-dimensional rules on the underlying sparse grid. Like TP method, the SG method is quite general and easy to implement. But unlike the TP method, its computational cost is much lower because the number of its required function evaluations grows exponentially with much smaller base.

The goal of this paper is to apply the MDI strategy to the SG method, the resulting algorithm, which is called the SG-MDI algorithm, provides a fast algorithm for an efficient implementation of the SG method. The SG-MDI method incorporates the MDI nested iteration idea into the sparse grid method which allows the reuse of the computation in the function evaluations at the integration points as much as possible. This saving significantly reduces the overall computational cost for implementing the sparse grid method from an exponential growth to a low-order polynomial growth.

The rest of this paper is organized as follows. In section 2, we first briefly recall the formulation of the sparse grid method and then introduce our SG-MDI algorithm in two and three dimensions to explain the main ideas of the algorithm as well as generalize it to arbitrary dimensions. In section 4, we present various numerical experiments to test the performance of the proposed SG-MDI algorithm and compare its performances to the standard SG method and the classical MC method. These numerical tests show that the SG-MDI algorithm is much faster in low and medium dimensions (i.e. $d \lesssim 100$) and in very high dimensions (i.e. $d \approx 1000$). It still works even when the MC method fails to compute. In section 5, we provide numerical tests to measure the influence of parameters in the proposed SG-MDI algorithm, including dependencies on the choice of underlying 1-d quadrature rule, the precision layer, and the choice of iteration step size. In section 6 we numerically estimate the computational complexity of the SG-MDI algorithm. This is done by using standard regression technique to discover the relationship between CPU time and dimension. It is showed that the CPU time grows at most in a polynomial order $O(d^3N)$, where $d$ and $N$ stand for respectively the dimension of integration domain and the maximum number of integration points in each coordinate direction. As a result, the SG-MDI algorithm can easily compute intermediate and very high-dimensional integrals on common desktop computers. Finally, the paper is concluded with some concluding remarks given in section 7.

**2. Preliminaries.** In this paper, $f(\mathbf{x}) : \bar{\Omega} \to R$ denote a generic continuous function on $\bar{\Omega} \subset \mathbb{R}^d$ for $d >> 1$, then $f(\mathbf{x})$ has pointwise values at every $\mathbf{x} = (x_1, x_2, \cdots, x_d) \in \bar{\Omega}$. Without loss of the generality, we set $\Omega := [-1, 1]^d$ and consider the basic and fundamental problem of computing the integral

$$(2.1) \qquad I_d(f) := \int_\Omega f(\mathbf{x}) d\mathbf{x}.$$

As mentioned in Section 1, the goal of this paper is to develop a fast algorithm for computing above integral based on the sparse grid methodology. To that end, below we briefly recall the necessary elements of sparse grid methods.

**2.1. The sparse grid method.** We now recall the formulation of the sparse grid method for approximating (2.1) and its tensor product reformulation formula which will be crucially used later in the formulation of our fast SG-MDI algorithm.

For each positive integer index $l \geq 1$, let $n_l$ be a positive integer which denotes the number of grid points at level $l$ and

$$(2.2) \qquad \Gamma_1^l := \left\{ x_1^l < x_2^l < \cdots < x_{n_l}^l \right\} \subset [-1, 1]$$

denote a sequence of the level $l$ grid points in $[-1, 1]$. The grid sets $\{\Gamma_1^l\}$ is said to be nested provided that $\Gamma_1^l \subset \Gamma_1^{l+1}$. The best known example of the nested grids is the following dyadic grids:

$$(2.3) \qquad \Gamma_1^l = \left\{ x_i^l := \frac{i}{2^{l-1}} : i = -2^{l-1}, \cdots, 0, 1, \cdots, 2^{l-1} \right\}.$$

For a given positive integer $q$, the tensor product $G_d^q := \Gamma_1^q \times \Gamma_1^q \times \cdots \times \Gamma_1^q$ then yields a standard tensor product grid on $\Omega = [-1, 1]^d$. Notice that here the $q$th level is used in each coordinate direction. To reduce the number of grid points in $G_d^q$, the sparse grid idea is to restrict the total level to be $q$ in the sense that $q = l_1 + l_2 + \cdots + l_d$, where $l_i$ is the level used in the $i$th coordinate direction, its corresponding tensor product grid is $\Gamma_1^{l_1} \times \Gamma_1^{l_2} \times \cdots \times \Gamma_1^{l_d}$. Obviously, the decomposition $q = l_1 + l_2 + \cdots + l_d$ is not unique, so all such decomposition must be considered. The union

$$(2.4) \qquad \Gamma_d^q := \bigcup_{l_1 + \cdots + l_d = q} \Gamma_1^{l_1} \times \cdots \times \Gamma_1^{l_d}$$

then yields the famous Smolyak sparse grid on $\Omega = [-1, 1]^d$ (cf. [7]). We remark that the underlying idea of going from $G_d^q$ to $\Gamma_d^q$ is exactly same as going from $Q_q(\Omega)$ to $P_q(\Omega)$, where $Q_q$ denotes the set of polynomials whose degrees in all coordinate directions not exceeding $q$ and $P_q$ denotes the set of polynomials whose total degrees not exceeding $q$.

After having introduced the concept of sparse grids, we then can define the sparse grid quadrature rule. For a univariate function $g$ on $[-1, 1]$, we consider $d$ one-dimensional quadrature formula

$$(2.5) \qquad \mathcal{J}_1^{l_i}(g) := \sum_{j=1}^{n_{l_i}} w_j^{l_i} g(x_j^{l_i}), \qquad i = 1, \cdots, d,$$

where $\{x_j^{l_i}, j = 1, \cdots, n_{l_i}\}$ and $\{w_j^{l_i}, j = 1, \cdots, n_{l_i}\}$ denote respectively the integration points/nodes and weights of the quadrature rule, and $n_{l_i}$ denotes the number of integration points in the $i$th coordinate direction in $[-1,1]$. Define

$$\mathbf{l} := [l_1, l_2, \cdots, l_d], \qquad |\mathbf{l}| := l_1 + l_2 + \cdots + l_d,$$

$$\mathcal{N}_{d,s} := \left\{ \mathbf{l} = [l_1, \cdots, l_d] : |\mathbf{l}| = s, s \geq d \right\}.$$

For example $\mathcal{N}_{2,4} = \left\{ [1,3], [2,2], [3,1] \right\}$.

Then, the sparse grid quadrature rule with accuracy level $q \in \mathbb{N}$ for $d$-dimensional integration (2.1) on $[-1,1]^d$ is defined as (cf. [11])

$$(2.6) \quad Q_d^q(f) := \sum_{q-d+1 \leq |\mathbf{l}| \leq q} (-1)^{q-|\mathbf{l}|} \left( \begin{array}{c} d-1 \\ q-|\mathbf{l}| \end{array} \right) \sum_{\mathbf{l} \in \mathcal{N}_{d,|\mathbf{l}|}} \left( \mathcal{J}_1^{l_1} \otimes \mathcal{J}_1^{l_2} \otimes \cdots \otimes \mathcal{J}_1^{l_d} \right) f,$$

where

$$(2.7) \qquad \left( \mathcal{J}_1^{l_1} \otimes \mathcal{J}_1^{l_2} \otimes \cdots \otimes \mathcal{J}_1^{l_d} \right) f := \sum_{j_1=1}^{n_{l_1}} \cdots \sum_{j_d=1}^{n_{l_d}} w_{j_1}^{l_1} \cdots w_{j_d}^{l_d} f\left( x_{j_1}^{l_1}, \cdots, x_{j_d}^{l_d} \right).$$

We note that each term $\left( \mathcal{J}_1^{l_1} \otimes \mathcal{J}_1^{l_2} \otimes \cdots \otimes \mathcal{J}_1^{l_d} \right) f$ in (2.7) is the tensor product quadrature rule which uses $n_{l_i}$ integration points in $i$th coordinate direction. To write $Q_d^q(f)$ more compactly, we set

$$n_d^q := \sum_{q-d+1 \leq |\mathbf{l}| \leq q} n_{l_1} \ldots n_{l_d},$$

which denotes the total number of integration points in $\Omega = [-1,1]^d$. Let $w_k, k = 1, \cdots, n_d^q$ denote the corresponding weights and define the bijective mapping

$$\left\{ \mathbf{x}^k : k = 1, \cdots, n_d^q \right\} \longrightarrow \left\{ (x_{j_1}^{l_1}, \cdots, x_{j_d}^{l_d}) : j_i = 1, \cdots, n_{l_i}, q-d+1 \leq |\mathbf{l}| \leq q \right\}.$$

Then, the sparse grid quadrature rule $Q_d^q(f)$ can be rewritten as

$$(2.8) \qquad\qquad Q_d^q(f) = \sum_{k=1}^{n_d^q} w_k f(\mathbf{x}^k).$$

We also note that some weights $w_k$ may become negative even though the one-dimensional weights $w_{j_1}^{l_1}, \cdots, w_{j_d}^{l_d}$ are positive. Therefore, it is no longer possible to interpret $Q_d^q(f)$ as a discrete probability measure. Moreover, the existence of negative weights in (2.8) may cause numerical cancellation, hence, loss of significant digits. To circumvent such a potential cancellation, it is recommended in [10] that the summation is carried out by coordinates, this then leads to the following tensor product reformulation of $Q_d^q(f)$:

$$(2.9) \qquad Q_d^q(f) = \sum_{l_1=1}^{q-1} \sum_{l_2=1}^{\gamma_1^q} \cdots \sum_{l_d=1}^{\gamma_{d-1}^q} \sum_{j_1=1}^{n_{l_1}} \cdots \sum_{j_d=1}^{n_{l_d}} w_{j_1}^{l_1} \cdots w_{j_d}^{l_d} f(x_{j_1}^{l_1}, \cdots, x_{j_d}^{l_d}),$$

where the upper limits are defined recursively as

$$\gamma_0^q := q, \quad \gamma_j^q := \gamma_{j-1}^q - l_j, \quad j = 1, 2, \cdots, d.$$

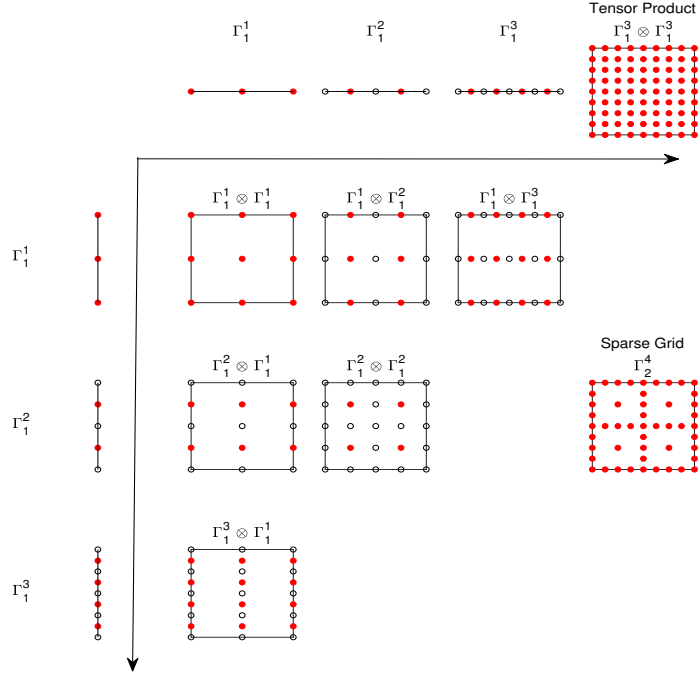

Fig. 2.1: Construction of a nested sparse grid in two dimensions.

In the nested mesh case (i.e., $\Gamma_1^l$ is nested), non-nested integration points are selected to form the sparse grid. We remark that different 1-d quadrature rules in (2.5) will lead to different sparse grid methods in (2.6). We also note that the tensor product reformulation (2.9) will play a crucial role later in the construction of our SG-MDI algorithm.

Figure 2.1 demonstrates the construction of a sparse grid according to the Smolyak rule when $d = 2$ and $q = 4$. The meshes are nested, namely, $\Gamma_1^{l_1} \subset \Gamma_1^{l_{i+1}}$. The 1-d integration-point sequence $\Gamma_1^{l_1}$ ($l_1 = 1, 2, 3$, and $n_{l_1} = 3, 5, 9$) and $\Gamma_1^{l_2}$ ($l_2 = 1, 2, 3$, and $n_{l_2} = 3, 5, 9$) are shown at the top and left of the figure, and the tensor product points $\Gamma_1^3 \otimes \Gamma_1^3$ are shown in the upper right corner. From (2.6) we see that the sparse grid rule is a combination of the low-order tensor product rule on $\Gamma_1^i \otimes \Gamma_1^j$ with $3 \leq i + j \leq 4$. The point sets of these products and the resulting sparse grid $\Gamma_2^4$ are shown in the lower half of the figure. We notice that some points in the sparse grid $\Gamma_2^4$ are repeatedly used in $\Gamma_1^i \otimes \Gamma_1^j$ with $3 \leq i + j \leq 4$. Consequently, we would avoid the repeated points (i.e., only using the red points in Figure 2.1) and use the reformulation (2.9) which does not involve repetition in the summation.

**2.2. Examples of sparse grid methods.** As mentioned earlier, different 1-d quadrature rules in (2.5) lead to different sparse grid quadrature rules in (2.9). Below we introduce four widely used sparse grid quadrature rules which will be the focus of this paper.

**Example 1: The classical trapezoidal rule.** The 1-d trapezoidal rule is

defined by (cf. [3])

$$\mathcal{J}_1^q(g) = \sum_{l=1}^{q-1} \sum_{j=1}^{n_l} {}''\, w_j^l\, g(x_j^l) \tag{2.10}$$

with $q \geq 2$ and $n_1 = 1, n_l = 2^{l-1} + 1$, $w_j^l = 2^{2-l}$, $x_j^l = (j-1) \cdot 2^{2-l} - 1$. Where the $\sum{}''$ indicates that the first and last terms in the sum are halved. It is well-known that there holds the following error estimate: if $g \in C^2$, then

$$\left| I_1(g) - \mathcal{J}_1^q(g) \right| \leq C\, 2^{-2q}.$$

**Example 2: The classical Clenshaw-Curtis rule.** This quadrature rule reads as follows (cf. [19, 20]):

$$\mathcal{J}_1^q(g) = \sum_{l=1}^{q-1} \sum_{j=1}^{n_l} w_j^l\, g(x_j^l) \tag{2.11}$$

with $q \geq 2$, $n_1 = 1, n_l = 2^{l-1} + 1$, $x_j^l = -\cos(\pi(j-1)/(n_l - 1))$, $j = 1, \cdots, n_l$ and the weights

$$w_1^l = w_{n_l}^l = \frac{1}{n_l(n_l - 2)}, \quad w_j^l = \frac{2}{n_l - 1}\left(1 + 2\sum_{k=1}^{\frac{n_l-1}{2}}{}'\,\frac{1}{1 - 4k^2}\cos\frac{2\pi(j-1)k}{n_l - 1}\right)$$

for $2 \leq j \leq n_l - 1$. Where $\sum{}'$ indicates that the last term in the summation is halved. It is well-known that there holds the following error estimate: if $g \in C^r$, then

$$\left| I_1(g) - \mathcal{J}_1^q(g) \right| \leq C\, (n_1^q)^{-r}, \qquad \text{where} \quad n_1^q := n_1 + n_2 + \cdots + n_{q-1}.$$

**Example 3: The Gauss-Patterson rule.** This quadrature rule is defined by (cf. [18, 23])

$$\mathcal{J}_1^q(g) = \sum_{l=1}^{q-1} \sum_{j=1}^{n_l} w_j^l\, g(x_j^l) \tag{2.12}$$

with $q \geq 2$, $n_l = 2^l(n+1) - 1$, and $\{x_j^l\}_{j=1}^{n_l}$ being the union of the zeroes of the polynomial $P_n(x)$ and $G_i(x), 1 \leq i < l$, where $P_n(x)$ is the $n$-th order Legendre polynomial and $G_1(x)$ is the $(n+1)$-th order Stieltjes polynomial and the $G_i(x)$ is orthogonal to all polynomials of degree less than $2^{l-1}(n+1)$ with respect to the weight function $P_n(x)(\prod_{j=1}^{i-1} G_j(x))$. $\{x_j^l\}_{j=1}^{n_l}$ are defined similarly to the Gauss-Legendre case, see [23] for details. Gauss-Patterson rules are a sequence of nested quadrature formulas with the maximal order of exactness. It is well-known that there holds the following error estimate: if $g \in C^r$, then

$$\left| I_1(g) - \mathcal{J}_1^q(g) \right| \leq C\, (n_1^q)^{-r}, \qquad \text{where} \quad n_1^q := n_1 + n_2 + \cdots + n_{q-1}.$$

**Example 4: The classical Gauss-Legendre rule.** This classical quadrature rule is defined by

$$\mathcal{J}_1^q(g) = \sum_{l=1}^{q-1} \sum_{j=1}^{n_l} w_j^l\, g(x_j^l) \tag{2.13}$$

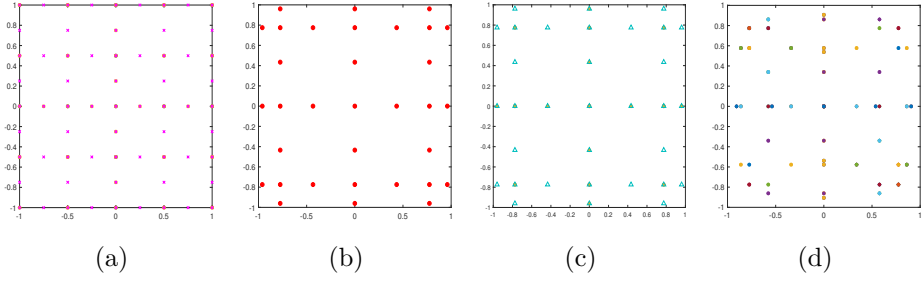

Fig. 2.2: Sparse grids corresponding to the trapezoidal rule (a), Clenshaw-Curtis rule (b), Gauss-Patterson rule (c), and Gauss-Legendre rule (d) when $d = 2, q = 6$.

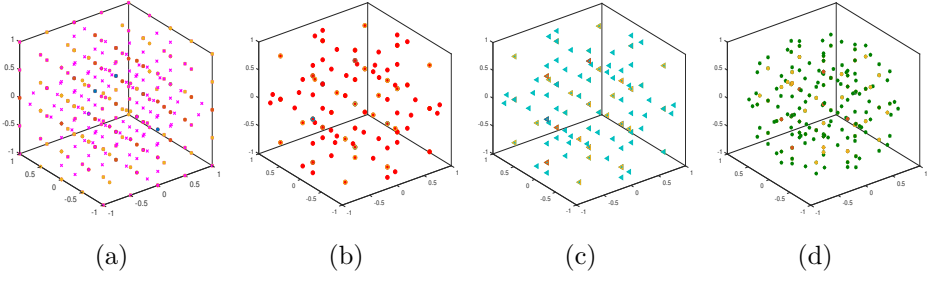

Fig. 2.3: Sparse grids corresponding to the trapezoidal rule (a), Clenshaw-Curtis rule (b), Gauss-Patterson rule (c), and Gauss-Legendre rule (d) when $d = 3, q = 7$.

with $q \geq 2, n_l = l$ for $l \geq 1$. $\{x_j^l\}_{j=1}^{n_l}$ are the zeroes of the $n_l$th order Legendre polynomial $P_{n_l}(x)$, and $\{w_j^l\}_{j=1}^{n_l}$ are the corresponding weights. It is well-known that there holds the following error estimate: if $g \in C^r$, then

$$\left| I_1(g) - \mathcal{J}_1^q(g) \right| \leq C \left( n_1^q \right)^{-r}, \qquad \text{where} \quad n_1^q := n_1 + n_2 + \cdots + n_{q-1}.$$

Figures 2.2 and 2.3 show the resulting sparse grids of the above four examples in both 2-d and 3-d cases with $q = 6$ and $q = 7$ respectively. We note that these four sparse grids have different structures.

We conclude this section by remarking that the error estimates of the above quadrature rules can be easily translate to error estimates for the sparse grid method (2.9). For example, in the case of the Clenshaw-Curtis, Gauss-Patterson, and Gauss-Legendre quadrature rule, there holds (cf. [4] )

$$(2.14) \qquad \left| I_d(f) - Q_d^q(f) \right| \leq C(n_d^q)^{-r} \cdot (\log(n_d^q))^{(d-1)(r+1)} \qquad \forall f \in W_d^r,$$

where

$$n_d^q := \sum_{q-d+1 \leq |\mathbf{l}| \leq q} n_{l_1} \ldots n_{l_d},$$

$$W_d^r := \left\{ f : \Omega \to \mathbb{R}; \ \left| \frac{\partial^{|\mathbf{l}|} f}{\partial x_1^{l_1} \cdots \partial x_d^{l_d}} \right|_\infty < \infty, \ l_i \leq r, i = 1, 2, \cdots, d \right\}.$$

We note that the above estimate indicates that the error of the sparse grid method still deteriorates exponentially in the dimension $d$, but with a much smaller base $\log(n_d^q)$.

**3. The SG-MDI algorithm.** The goal of this section is to present an efficient and fast implementation algorithm (or solver), called the SG-MDI algorithm, for evaluating the sparse grid quadrature rule (2.6) via its reformulation (2.9) in order to circumvent the curse of dimensionality which hampers the usage of the sparse grid method (2.6) in high dimensions. To better understand the main idea of the SG-MDI method, we first consider the simple two and three dimensional cases and then to formulate the algorithm in arbitrary dimensions.

**3.1. Formulation of the SG-MDI algorithm in two dimensions.** Let $d = 2$, $\Omega = [-1, 1]^2$, and $\mathbf{x} = (x_1, x_2) \in \Omega$. By Fubini's Theorem we have

$$(3.1) \qquad I_2(f) := \int_{[-1,1]^2} f(\mathbf{x}) \, d\mathbf{x} = \int_{-1}^1 \left( \int_{-1}^1 f(x_1, x_2) \, dx_1 \right) dx_2.$$

Then, the two-dimensional SG quadrature rule (2.9) takes the form

$$(3.2) \qquad Q_2^q(f) = \sum_{l_1=1}^{q-1} \sum_{l_2=1}^{\gamma_1^q} \sum_{j_1=1}^{n_{l_1}} \sum_{j_2=1}^{n_{l_2}} w_{j_1}^{l_1} w_{j_2}^{l_2} f(x_{j_1}^{l_1}, x_{j_2}^{l_2}),$$

where $\gamma_1^q = q - l_1$. Motivated by (and mimicking) the Fubini's formula (3.1), we rewrite (3.2) as

$$(3.3) \qquad \begin{aligned} Q_2^q(f) &= \sum_{l_2=1}^{\gamma_1^q} \sum_{j_2=1}^{n_{l_2}} w_{j_2}^{l_2} \left( \sum_{l_1=1}^{q-1} \sum_{j_1=1}^{n_{l_1}} w_{j_1}^{l_1} f(x_{j_1}^{l_1}, x_{j_2}^{l_2}) \right) \\ &= \sum_{l_2=1}^{\gamma_1^q} \sum_{j_2=1}^{n_{l_2}} w_{j_2}^{l_2} f_1(x_{j_2}^{l_2}), \end{aligned}$$

where

$$f_1(s) := \sum_{l_1=1}^{q-1} \sum_{i=1}^{n_{l_1}} w_i^{l_1} f(x_i^{l_1}, s).$$

We note that the evaluation of $f_1(x_2)$ is amount to applying the 1-d formula (2.5) to approximate the integral $\int_{-1}^1 f(x_1, x_2) dx_1$. However, the values of $\{f_1(x_{j_2}^{l_2})\}$ will not be computed by the 1-d quadrature rule in our SG-MDI algorithm, instead, $f_1$ is formed as a symbolic function, so the 1-d quadrature rule can be called on $f_1$. Therefore, we still use the SG method to select the integration points, and then use our SG-MDI algorithm to perform function evaluations at the integration points collectively to save computation, which is the main idea of the SG-MDI method.

Let $W$ and $X$ denote the weight and node vectors of the 1-d quadrature rule on $[-1, 1]$, $n_{l_i}$ represents the number of integration points in the $x_i$ direction and we use a parameter $r \in \{1, 2, 3, 4\}$ to indicate one of the four quadrature rule. The following algorithm implements the sparse grid quadrature formula (3.3).

---

**Algorithm 3.1** 2d-SG-MDI($f$, $\Omega$, $r$, $q$, $X$, $W$)

---

Initialize $Q = 0$, $f_1 = 0$.
**for** $l_1 = 1 : q$ **do**
$\quad$ **for** $j_1 = 1 : n_{l_1}$ **do**
$\quad\quad$ $f_1 = f_1 + W_{j_1}^{l_1} f((X_{j_1}^{l_1}, \cdot))$.
$\quad$ **end for**
**end for**
**for** $l_2 = 1 : q - l_1$ **do**
$\quad$ **for** $j_2 = 1 : n_{l_2}$ **do**
$\quad\quad$ $Q = Q + W_{j_2}^{l_2} \cdot f_1(X_{j_2}^{l_2})$.
$\quad$ **end for**
**end for**
**return** $Q$.

---

We note that the first do-loop forms the symbolic function $f_1$ which encodes all computations involving the $x_1$-components at all integration points. The second do-loop evaluates the 1-d quadrature rule for the function $f_1$. As mentioned above, in this paper we only focus on the four well-known 1-d quadrature rules: (i) trapezoidal rule; (ii) Clenshaw-Curtis rule; (iii) Gauss-Patterson rule; (iv) Gauss-Legendre rule. They will be represented respectively by $r = 1, 2, 3, 4$.

**3.2. Formulation of the SG-MDI algorithm in three dimensions.** In the subsection we extend the formulation of the above 2d-SG-MDI algorithm to the 3-d case by highlighting its main steps, in particular, how the above 2-d algorithm is utilized. First, recall that Fubini's Theorem is given by

$$(3.4) \qquad I_3(f) := \int_{[-1,1]^3} f(\mathbf{x}) \, d\mathbf{x} = \int_{[-1,1]^2} \left( \int_{-1}^{1} f(\mathbf{x}) \, dx_1 \right) d\mathbf{x}',$$

where $\mathbf{x}' = (x_2, x_3)$.

Second, notice that the SG quadrature rule (2.9) in 3-d takes the form

$$(3.5) \qquad Q_3^q(f) = \sum_{l_1=1}^{q-1} \sum_{l_2=1}^{\gamma_1^q} \sum_{l_3=1}^{\gamma_2^q} \sum_{j_1=1}^{n_{l_1}} \sum_{j_2=1}^{n_{l_2}} \sum_{j_3=1}^{n_{l_3}} w_{j_1}^{l_1} w_{j_2}^{l_2} w_{j_3}^{l_3} \, f(x_{j_1}^{l_1}, x_{j_2}^{l_2}, x_{j_3}^{l_3}).$$

where $\gamma_2^q = q - l_1 - l_2$. Mimicking the above Fubini's formula, we rewrite (3.5) as

$$(3.6) \qquad Q_3^q(f) = \sum_{l_3=1}^{\gamma_2^q} \sum_{j_3=1}^{n_{l_3}} \sum_{l_2=1}^{\gamma_1^q} \sum_{j_2=1}^{n_{l_2}} w_{j_3}^{l_3} w_{j_2}^{l_2} \left( \sum_{l_1=1}^{q-1} \sum_{j_1=1}^{n_{l_1}} w_{j_1}^{l_1} \, f(x_{j_1}^{l_1}, x_{j_2}^{l_2}, x_{j_3}^{l_3}) \right)$$
$$= \sum_{l_3=1}^{\gamma_2^q} \sum_{j_3=1}^{n_{l_3}} \sum_{l_2=1}^{\gamma_1^q} \sum_{j_2=1}^{n_{l_2}} w_{j_3}^{l_3} w_{j_2}^{l_2} \, f_2(x_{j_2}^{l_2}, x_{j_3}^{l_3}),$$

where

$$(3.7) \qquad f_2(s, t) := \sum_{l_1=1}^{q-1} \sum_{j_1=1}^{n_{l_1}} w_{j_1}^{l_1} f(x_{j_1}^{l_1}, s, t).$$

We note that $f_2$ is formed as a symbolic function in our SG-MDI algorithm and the right-hand side of (3.6) is viewed as a 2-d sparse grid quadrature formula for $f_2$, it can be computed either directly or recursively by using **Algorithm 3.1**. The following algorithm implements the SG quadrature formula (3.6).

---

**Algorithm 3.2** 3d-SG-MDI($f$, $\Omega$, $r$, $q$, $X$, $W$)

---

    Initialize $Q = 0, f_2 = 0$.
    **for** $l_1 = 1 : q$ **do**
        **for** $j_1 = 1 : n_{l_1}$ **do**
            $f_2 = f_2 + W_{j_1}^{l_1} \cdot f((X_{j_1}^{l_1}, \cdot, \cdot))$.
        **end for**
    **end for**
    $\Omega_2 = P_3^2 \Omega, \gamma_1^q = q - l_1$.
    $Q =$ 2d-SG-MDI($f_2, \Omega_2, r, \gamma_1^q, X, W$).
    **return** $Q$.

---

Where $P_3^2$ denotes the orthogonal projection (or natural embedding): $x = (x_1, x_2, x_3)$ $\rightarrow \mathbf{x}' = (x_2, x_3)$, $W$ and $X$ stand for the weight and node vectors of the underlying 1-d quadrature rule.

From **Algorithm 2.2** we can see the mechanism of the SG-MDI algorithm. It is based on two main ideas: (i) to use the sparse grid approach to select integration points; (ii) to use the discrete Fubini formula to efficiently compute the total sparse grid sum by reducing it to calculation of a low-dimensional (i.e., 2-d) sparse grid sum, which allows us to recursively call the low-dimensional SG-MDI algorithm.

**3.3. Formulation of the SG-MDI algorithm in arbitrary d-dimensions.** The goal of this subsection is to extend the 2-d and 3-d SG-MDI algorithms to arbitrary d-dimensions. We again start with recalling the d-dimensional Fubini's Theorem

$$(3.8) \qquad I_d(f) = \int_\Omega f(\mathbf{x}) \, d\mathbf{x} = \int_{\Omega_{d-m}} \left( \int_{\Omega_m} f(\mathbf{x}) \, d\mathbf{x}'' \right) d\mathbf{x}',$$

where $1 \leq m < d$, $\Omega = [-1, 1]^d, \Omega_m = \rho_d^m \Omega = [-1, 1]^m$ and $\Omega_{d-m} = P_d^{d-m}\Omega = [-1, 1]^{d-m}$ in which $\rho_d^m$ and $P_d^{d-m}$ denote respectively the orthogonal projections (or natural embeddings): $\mathbf{x} = (x_1, x_2, \cdots, x_d) \rightarrow \mathbf{x}'' = (x_1, x_2, \cdots, x_m)$ and $\mathbf{x} = (x_1, x_2, \cdots, x_d) \rightarrow \mathbf{x}' = (x_{m+1}, x_{m+2}, \cdots, x_d)$. The integer $m$ denotes the dimension reduction step length in our algorithm.

Mimicking the above Fubini's Theorem, we rewrite the d-dimensional SG quadrature rule (2.9) as follows:

$$(3.9)$$
$$Q_d^q(f) = \sum_{l_d=1}^{\gamma_{d-1}^q} \sum_{j_d=1}^{n_{l_d}} \cdots \sum_{l_{m+1}=1}^{\gamma_m^q} \sum_{j_{m+1}=1}^{n_{l_{m+1}}} w_{j_{m+1}}^{l_{m+1}} \cdots w_{j_d}^{l_d} \left( \sum_{l_m=1}^{\gamma_{m-1}^q} \cdots \sum_{j_1=1}^{n_{l_1}} w_{j_1}^{l_1} \cdots w_{j_m}^{l_m} f(x_{j_1}^{l_1}, \cdots, x_{j_d}^{l_d}) \right)$$
$$= \sum_{l_d=1}^{\gamma_{d-1}^q} \sum_{j_d=1}^{n_{l_d}} \cdots \sum_{l_{m+1}=1}^{\gamma_m^q} \sum_{j_{m+1}=1}^{n_{l_{m+1}}} w_{j_{m+1}}^{l_{m+1}} \cdots w_{j_d}^{l_d} \, f_{d-m}(x_{j_{m+1}}^{l_{m+1}}, \cdots, x_{j_d}^{l_d}),$$

where

$$(3.10)$$
$$f_{d-m}(s_1, \cdots, s_{d-m}) = \sum_{l_m=1}^{\gamma_{m-1}^q} \sum_{j_m=1}^{n_{l_m}} \cdots \sum_{l_1=1}^{q-1} \sum_{j_1=1}^{n_{l_1}} w_{j_1}^{l_1} \cdots w_{j_m}^{l_m} \, f(x_{j_1}^{l_1}, \cdots, x_{j_m}^{l_m}, s_1, \cdots, s_{d-m}).$$

We note that in our SG-MDI algorithm $f_{d-m}$ defined by (3.10) is a symbolic function and the right-hand side of (3.9) is a $(d-m)$-order multi-summation, which itself can be evaluated by employing the dimension reduction strategy. Dimensionality can be reduced by iterating $k := [\frac{d}{m}]$ times until $d - km \leq m$. To implement this process, we introduce the following conventions.

- If $t = 1$, set SG-MDI$(t, f_t, \Omega_t, m, s, r, q, X, W) := \mathcal{J}_1^q(f)$, which is computed by using the one-dimensional quadrature rule (2.5).
- If $t = 2$, set SG-MDI$(t, f_t, \Omega_t, m, s, r, q, X, W) :=$ 2d-SG-MDI$(f_t, \Omega_t, r, q, X, W)$.
- If $t = 3$, set SG-MDI$(t, f_t, \Omega_t, m, s, r, q, X, W) :=$ 3d-SG-MDI$(f_t, \Omega_t, r, q, X, W)$.

We note that when $t = 1, 2, 3$, the parameter $m$ becomes a dummy variable and can be given any value. Recall that $P_t^{t-m}$ denotes the natural embedding from $\mathbb{R}^t$ to $\mathbb{R}^{t-m}$ by deleting the first $m$ components of vectors in $\mathbb{R}^t$. The following algorithm implements the sparse grid quadrature via (3.9).

---

**Algorithm 3.3** SG-MDI$(d, f, \Omega, m, s, r, q, X, W)$

---

$\Omega_d = \Omega$, $f_d = f$, $k = [\frac{d}{m}]$, $\gamma_d^q = q - 1$.
**for** $t = d : -m : d - km$ (the index is decreased by $m$ at each iteration) **do**
    $\Omega_{d-m} = P_t^{t-m} \Omega_t$,     $\gamma_{d-m}^q = \gamma_t^q - l_1 - \cdots - l_m$.
    (Construct symbolic function $f_{t-m}$ by (3.11) below).
    SG-MDI$(t, f_t, \Omega_t, m, s, r, \gamma_t^q, X, W) :=$ SG-MDI$(t-m, f_{t-m}, \Omega_{t-m}, m, s, r, \gamma_{t-m}^q, X, W)$
**end for**
$d = d - km$, $s = 1$(or 2, 3), $f_d = f_t$, $k_1 = [\frac{d-km}{s}]$.
**for** $t = d : s : d - k_1 s$ (the index is decreased by $s$ at each iteration) **do**
    $\Omega_{d-s} = P_t^{t-s} \Omega_t$,     $\gamma_{d-s}^q = \gamma_t^q - l_1 - \cdots - l_s$.
    (Construct symbolic function $f_{t-s}$ by (3.11) below).
    SG-MDI$(t, f_t, \Omega_t, m, s, r, \gamma_t^q, X, W) :=$ SG-MDI$(t - s, f_{t-s}, \Omega_{t-s}, m, s, r, \gamma_{t-s}^q, X, W)$
**end for**
$Q =$ SG-MDI$(d - k_1 s, f_{d-k_1 s}, \Omega_{d-k_1 s}, m, s, r, \gamma_{d-k_1 s}^q, X, W)$.
**return** $Q$.

---

Where

$$(3.11)$$
$$f_{t-m}(s_1, \cdots, s_{t-m}) = \sum_{l_m=1}^{\gamma_{m-1}^q} \sum_{j_m=1}^{n_{l_m}} \cdots \sum_{l_1=1}^{q-1} \sum_{j_1=1}^{n_{l_1}} w_{j_1}^{l_1} \cdots w_{j_m}^{l_m} \cdot f(x_{j_1}^{l_1}, \cdots, x_{j_m}^{l_m}, s_1, \cdots, s_{t-m}).$$

We remark that **Algorithm 2.3** recursively generates a sequence of symbolic functions $f_d, \cdots, f_{d-km}, f_{d-km-s}, \cdots, f_{d-km-k_1 s}$, each function has $m$ fewer arguments than its predecessor. As mentioned earlier, our SG-MDI algorithm does not perform the function evaluations at all integration points independently, but rather iteratively along $m$-coordinate directions, hence, the function evaluation at any integration point is not completed until the last step of the algorithm is executed. As a result, many computations are reused in each iteration, which is the main reason for the computation saving and to achieve a faster algorithm.

**4. Numerical performance tests.** In this section, we present extensive numerical tests to guage the performance of the proposed SG-MDI algorithm and to compare it with the standard sparse grid (SG) and classical Monte Carlo (MC) methods. All numerical tests show that SG-MDI algorithm outperforms both SG and MC methods in low and medium high dimensions (i.e. $d \lesssim 100$), and can compute very high-dimensional (i.e. $d \approx 1000$) integrals while others fail.

All our numerical experiments are done in Matlab 9.4.0.813654(R2018a) on a desktop PC with Intel(R) Xeon(R) Gold 6226R CPU 2.90GHz and 32GB RAM.

**4.1. Two and three-dimensional tests.** We first test our SG-MDI algorithm on simple 2-d and 3-d examples and compare its performance (in terms of CPU time) with SG methods.

**Test 1.** Let $\Omega = [-1, 1]^2$ and consider the following 2-d integrands:

$$(4.1) \qquad f(x) := \exp\big(5x_1^2 + 5x_2^2\big); \qquad \widehat{f}(x) := \sin\big(2\pi + 10x_1^2 + 5x_2^2\big).$$

Let $q$ denote the accuracy level of the sparse grid. The larger $q$ is, the more integration points we need for the 1-d quadrature rule, and the higher the accuracy of the SG-MDI quadrature. The base 1-d quadrature rule is chosen to be the Gauss-Patterson rule, hence, parameter $r = 3$ in the algorithm.

Table 4.1 and 4.2 present the computational results (errors and CPU times) of the SG-MDI and SG methods for approximating $I_2(f)$ and $I_2(\widehat{f})$, respectively.

| Accuracy level ($q$) | Total nodes | SG-MDI | | SG | |
|---|---|---|---|---|---|
| | | Relative error | CPU time(s) | Relative error | CPU time(s) |
| 6 | 33 | $1.0349 \times 10^{-1}$ | 0.0512817 | $1.0349 \times 10^{-1}$ | 0.0078077 |
| 7 | 65 | $2.3503 \times 10^{-3}$ | 0.0623538 | $2.3503 \times 10^{-3}$ | 0.0084645 |
| 9 | 97 | $8.1019 \times 10^{-4}$ | 0.0644339 | $8.1019 \times 10^{-4}$ | 0.0095105 |
| 10 | 161 | $1.8229 \times 10^{-6}$ | 0.0724491 | $1.8229 \times 10^{-6}$ | 0.0106986 |
| 13 | 257 | $2.0720 \times 10^{-7}$ | 0.0913161 | $2.0720 \times 10^{-7}$ | 0.0135131 |
| 14 | 321 | $4.3279 \times 10^{-7}$ | 0.1072016 | $4.3279 \times 10^{-7}$ | 0.0155733 |

Table 4.1: Relative errors and CPU times of SG-MDI and SG tests with $m = 1$ for approximating $I_2(f)$.

| Accuracy level ($q$) | Total nodes | SG-MDI | | SG | |
|---|---|---|---|---|---|
| | | Relative error | CPU time(s) | Relative error | CPU time(s) |
| 9 | 97 | $4.7425 \times 10^{-1}$ | 0.0767906 | $4.7425 \times 10^{-1}$ | 0.0098862 |
| 10 | 161 | $1.4459 \times 10^{-3}$ | 0.0901238 | $1.4459 \times 10^{-3}$ | 0.0102700 |
| 13 | 257 | $1.9041 \times 10^{-5}$ | 0.1025934 | $1.9041 \times 10^{-5}$ | 0.0152676 |
| 14 | 321 | $2.3077 \times 10^{-5}$ | 0.1186194 | $2.3077 \times 10^{-5}$ | 0.0144737 |
| 16 | 449 | $3.1236 \times 10^{-6}$ | 0.1353691 | $3.1236 \times 10^{-6}$ | 0.0177445 |
| 20 | 705 | $2.4487 \times 10^{-6}$ | 0.1880289 | $2.4487 \times 10^{-6}$ | 0.0355606 |

Table 4.2: Relative errors and CPU times of SG-MDI and SG tests with $m = 1$ for approximating $I_2(\widehat{f})$.

From Table 4.1 and 4.2, we observe that these two methods use very little CPU time, although the SG method in both tests outperforms, but the difference is almost negligible, so both methods do well in 2-d case.

**Test 2.** Let $\Omega = [-1, 1]^3$ and we consider the following 3-d integrands:

$$(4.2) \qquad f(x) = \exp\left(5x_1^2 + 5x_2^2 + 5x_3^2\right), \qquad \widehat{f}(x) = \sin\left(2\pi + 10x_1^2 + 5x_2^2 + 5x_3\right).$$

We compute the integral of these two functions over $\Omega = [-1, 1]^3$ using the SG-MDI and SG methods. Likewise, let $q$ denote the accuracy level of the sparse grid, choose parameters $r = 3$ and $m = 1$ in the algorithm.

The test results are given in Table 4.3 and 4.4. We again observe that both methods use very little CPU time although the SG method again slightly outperforms in both tests. However, as $q$ increases, the number of integration points increases, and the CPU times used by these two methods get closer. We would like to point out that both methods are very efficient and their difference is negligible in the 3-d case.

| Accuracy level $(q)$ | Total nodes | SG-MDI | | SG | |
|---|---|---|---|---|---|
| | | Relative error | CPU time(s) | Relative error | CPU time(s) |
| 9 | 495 | $3.2467 \times 10^{-2}$ | 0.0669318 | $3.2467 \times 10^{-2}$ | 0.0235407 |
| 10 | 751 | $1.8956 \times 10^{-3}$ | 0.0886774 | $1.8956 \times 10^{-3}$ | 0.0411750 |
| 11 | 1135 | $3.9146 \times 10^{-4}$ | 0.0902602 | $3.9146 \times 10^{-4}$ | 0.0672375 |
| 13 | 1759 | $4.7942 \times 10^{-6}$ | 0.1088353 | $4.7942 \times 10^{-6}$ | 0.0589584 |
| 14 | 2335 | $1.8013 \times 10^{-6}$ | 0.1381728 | $1.8013 \times 10^{-6}$ | 0.0704032 |
| 15 | 2527 | $1.2086 \times 10^{-6}$ | 0.1484829 | $1.2086 \times 10^{-6}$ | 0.0902680 |
| 16 | 3679 | $3.6938 \times 10^{-7}$ | 0.1525743 | $3.6938 \times 10^{-7}$ | 0.1143728 |

Table 4.3: Relative errors and CPU times of SG-MDI and SG tests with $m = 1$ for approximating $I_3(f)$.

| Accuracy level $(q)$ | Total nodes | SG-MDI | | SG | |
|---|---|---|---|---|---|
| | | Relative error | CPU time(s) | Relative error | CPU time(s) |
| 12 | 1135 | $5.5057 \times 10^{-1}$ | 0.0921728 | $5.5057 \times 10^{-1}$ | 0.0495310 |
| 13 | 1759 | $8.9519 \times 10^{-3}$ | 0.1031632 | $8.9519 \times 10^{-3}$ | 0.0644124 |
| 15 | 2527 | $1.8063 \times 10^{-3}$ | 0.1771094 | $1.8063 \times 10^{-3}$ | 0.0891040 |
| 16 | 3679 | $1.1654 \times 10^{-4}$ | 0.1957219 | $1.1654 \times 10^{-4}$ | 0.1159222 |
| 17 | 4447 | $2.4311 \times 10^{-5}$ | 0.2053174 | $2.4311 \times 10^{-5}$ | 0.1443184 |
| 19 | 6495 | $5.4849 \times 10^{-6}$ | 0.4801467 | $5.4849 \times 10^{-6}$ | 0.2259950 |
| 20 | 8031 | $1.5333 \times 10^{-6}$ | 0.6777698 | $1.5333 \times 10^{-6}$ | 0.2632516 |

Table 4.4: Relative errors and CPU times of SG-MDI and SG tests with $m = 1$ for approximating $I_3(\widehat{f})$.

**4.2. High-dimensional tests.** In this section we evaluate the performance of the SG-MDI method for $d >> 1$. First, we test and compare the performance of SG-MDI and SG methods in computing Gaussian integrals for dimensions $2 \leq d \leq 20$ because $d \approx 20$ is the highest dimension that the SG method is able to compute a result on our computer. We then provide a performance comparison (in terms of CPU time) of the SG-MDI and classical Monte Carlo (MC) methods in computing high-dimensional integrals.

**Test 3.** Let $\Omega = [-1, 1]^d$ for $2 \leq d \leq 20$ and consider the following Gaussian integrand:

$$(4.3) \qquad f(\mathbf{x}) = \frac{1}{\sqrt{2\pi}} \exp\left(-\frac{1}{2}|\mathbf{x}|^2\right),$$

where $|\mathbf{x}|$ stands for the Euclidean norm of the vector $\mathbf{x} \in \mathbb{R}^d$.

We compute the integral $I_d(f)$ by using the SG-MDI and SG methods, as done in **Test 1-2**. Both methods are based on the same 1-d Gauss-Patterson rule (i.e., the parameter $r = 3$). We also set $m = 1, s = 1$ in the SG-MDI method and use two accuracy levels $q = 10, 12$, respectively.

| Dimension $(d)$ | Total nodes | SG-MDI | | SG | |
|---|---|---|---|---|---|
| | | Relative error | CPU time(s) | Relative error | CPU time(s) |
| 2 | 161 | $1.0163 \times 10^{-8}$ | 0.0393572 | $1.0163 \times 10^{-8}$ | 0.0103062 |
| 4 | 2881 | $2.0310 \times 10^{-8}$ | 0.0807326 | $2.0310 \times 10^{-8}$ | 0.0993984 |
| 8 | 206465 | $1.3429 \times 10^{-7}$ | 0.1713308 | $1.3429 \times 10^{-7}$ | 6.7454179 |
| 10 | 1041185 | $1.6855 \times 10^{-6}$ | 0.2553576 | $1.6855 \times 10^{-6}$ | 86.816883 |
| 12 | 4286913 | $1.8074 \times 10^{-5}$ | 0.3452745 | $1.8074 \times 10^{-5}$ | 866.1886366 |
| 14 | 5036449 | $2.1338 \times 10^{-4}$ | 0.4625503 | $2.1338 \times 10^{-4}$ | 6167.3838002 |
| 15 | 12533167 | $7.1277 \times 10^{-4}$ | 0.5867914 | failed | failed |

Table 4.5: Relative errors and CPU times of SG-MDI and SG tests with $m = 1, s = 1$, and $q = 10$ for approximating $I_d(f)$.

| Dimension $(d)$ | Total nodes | SG-MDI | | SG | |
|---|---|---|---|---|---|
| | | Relative error | CPU time(s) | Relative error | CPU time(s) |
| 2 | 161 | $1.0198 \times 10^{-8}$ | 0.0418615 | $1.0198 \times 10^{-8}$ | 0.0191817 |
| 4 | 6465 | $2.0326 \times 10^{-8}$ | 0.0704915 | $2.0326 \times 10^{-8}$ | 0.2067346 |
| 6 | 93665 | $3.0487 \times 10^{-8}$ | 0.0963325 | $3.0487 \times 10^{-8}$ | 3.1216913 |
| 8 | 791169 | $4.0881 \times 10^{-8}$ | 0.2233707 | $4.0881 \times 10^{-8}$ | 41.3632962 |
| 10 | 5020449 | $4.0931 \times 10^{-8}$ | 0.3740873 | $4.0931 \times 10^{-8}$ | 821.6461622 |
| 12 | 25549761 | $1.1560 \times 10^{-6}$ | 0.8169479 | $1.1560 \times 10^{-6}$ | 11887.797686 |
| 13 | 29344150 | $5.2113 \times 10^{-6}$ | 1.2380811 | failed | failed |

Table 4.6: Relative errors and CPU times of SG-MDI and SG tests with $m = 1, s = 1$, and $q = 12$ for approximating $I_d(f)$.

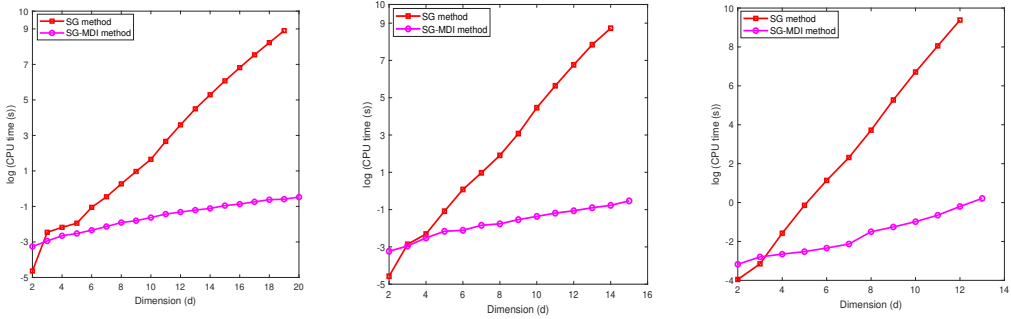

Fig. 4.1: CPU time comparison of SG-MDI and SG tests: $q = 8$ (left), $q = 10$ (middle), $q = 12$ (right).

Table 4.5 gives the relative error and CPU time for approximating $I_d(f)$ using SG-MDI and SG methods with accuracy level $q = 10$, and Table 4.6 gives the corresponding results for $q = 12$. We observe that the errors should be the same for both methods (since they use the same integration points), but their CPU times are quite different. The SG method is more efficient for $d \leq 4$ when $q = 8, 10$ and for $d \leq 3$ when $q = 12$, but the SG-MDI method excels for $d \geq 4$ and the winning margin becomes significant as $d$ and $q$ increase (also see Figure 4.1). For example, when $d = 14$ and $q = 10$, the CPU time required by the SG method is about 6167 seconds, which is about 2 hours, but the CPU time of the SG-MDI method is less than 1 second! Also, when $d = 13$ and $q = 12$, the SG method fails to compute the integral due to running out of computer memory because too large number of integration points must be saved and function evaluations must be performed, but the SG-MDI method only needs about 2 seconds to complete the computation!

The classical (and quasi) Monte Carlo (MC) method is often the preferred/default method for computing high-dimensional integrals. However, due to its low order of convergence, to achieve the accuracy, a large number of function evaluations are required at randomly sampled integration points and the number grows rapidly as dimension $d$ increases (due to the rapid growth of variance). Below we compare the performance of the SG-MDI (with parameters $r = 3, q = 10, m = 10, s = 1$) and classical MC method. In the test, when $d \geq 10$, we use the iteration step length $m > 1$ to iterate faster until $d - km \leq m$ to reach the stage 2 of the iteration. We refer the reader to Section 5.2 for a detailed analysis.

**Test 4.** Let $\Omega = [-1, 1]^d$ and choose the following integrands:

$$(4.4) \qquad f(x) = \prod_{i=0}^{d} \frac{1}{0.9^2 + (x_i - 0.6)^2}, \qquad \widehat{f}(x) = \frac{1}{\sqrt{2\pi}} \exp\left(-\frac{1}{2}|x|^2\right).$$

We use relative error as a criterion for comparison, that is, we determine a required (minimum) number of random sampling points for the MC method so that it produces a relative error comparable to that of the SG-MDI method. The computed results for $I_d(f)$ and $I_d(\widehat{f})$ are respectively given in Table 4.7 and 4.8.

From Table 4.7 and 4.8, we clearly see that there is a significant difference in the CPU time of these two methods for computing $I_d(f)$ and $I_d(\widehat{f})$. When $d > 30$, the classical MC method fails to produce a computed result with a relative error of order $10^{-5}$. As explained in [9], the MC method requires more than $10^{10}$ randomly sampled integration points and then needs independently to compute their function values, which is a tall order to do on a regular workstation.

| Dimension $(d)$ | MC | | SG-MDI | |
|---|---|---|---|---|
| | Relative error | CPU time(s) | Relative error | CPU time(s) |
| 5 | $1.3653 \times 10^{-5}$ | 62.1586394 | $1.3653 \times 10^{-5}$ | 0.0938295 |
| 10 | $2.0938 \times 10^{-5}$ | 514.1493073 | $2.0938 \times 10^{-5}$ | 0.1945813 |
| 20 | $4.2683 \times 10^{-5}$ | 1851.0461899 | $4.1876 \times 10^{-5}$ | 0.4204564 |
| 30 | $6.2814 \times 10^{-5}$ | 15346.222011 | $6.2814 \times 10^{-5}$ | 0.7692118 |
| 35 | $7.3283 \times 10^{-5}$ | failed | $7.3283 \times 10^{-5}$ | 0.9785784 |
| 40 | $8.3752 \times 10^{-5}$ | | $8.3752 \times 10^{-5}$ | 1.2452577 |
| 60 | $1.2562 \times 10^{-4}$ | | $1.2562 \times 10^{-4}$ | 2.5959174 |
| 80 | $1.6750 \times 10^{-4}$ | | $1.6750 \times 10^{-4}$ | 4.9092032 |
| 100 | $2.1235 \times 10^{-4}$ | | $2.1235 \times 10^{-4}$ | 8.1920274 |

Table 4.7: CPU times of the SG-MDI and MC tests with comparable relative errors for approximating $I_d(f)$.

| Dimension $(d)$ | MC | | SG-MDI | |
|---|---|---|---|---|
| | Relative error | CPU time(s) | Relative error | CPU time(s) |
| 5 | $9.4279 \times 10^{-7}$ | 85.2726354 | $9.4279 \times 10^{-7}$ | 0.0811157 |
| 10 | $1.6855 \times 10^{-6}$ | 978.1462121 | $1.6855 \times 10^{-6}$ | 0.295855 |
| 20 | $3.3711 \times 10^{-6}$ | 2038.138555 | $3.3711 \times 10^{-6}$ | 6.3939110 |
| 30 | $5.0567 \times 10^{-6}$ | 16872.143255 | $5.0567 \times 10^{-6}$ | 29.5098187 |
| 35 | $5.8995 \times 10^{-6}$ | failed | $5.8995 \times 10^{-6}$ | 62.0270714 |
| 40 | $6.7423 \times 10^{-6}$ | | $6.7423 \times 10^{-6}$ | 106.1616486 |
| 80 | $1.3484 \times 10^{-5}$ | | $1.3484 \times 10^{-5}$ | 1893.8402620 |
| 100 | $1.7825 \times 10^{-5}$ | | $1.7825 \times 10^{-5}$ | 3077.1890005 |

Table 4.8: CPU times of the SG-MDI and MC tests with comparable relative errors for approximating $I_d(\widehat{f})$.

Next, we come to address a natural question that asks how high the dimension $d$ can be handled by the SG-MDI method. Obviously, the answer must be computer-dependent, and that given below is obtained using the workstation at our disposal.

**Test 5.** Let $\Omega = [-1, 1]^d$ and consider the following integrands:

$$(4.5) \qquad f(x) = \frac{1}{2^d} \exp\Big(\sum_{i=1}^{d} (-1)^{i+1} x_i\Big), \qquad \widehat{f}(x) = \prod_{i=0}^{d} \frac{1}{0.9^2 + (x_i - 0.6)^2}.$$

We then compute $I_d(f)$ and $I_d(\widehat{f})$ using the SG-MDI algorithm for an increasing sequence of $d$ up to 1000 with parameters $r = 3, q = 10, m = 10, s = 1$. The computed results are shown in Tables 4.9. We stop the computation at $d = 1000$ since it is already quite high and use $q = 10, m = 10, s = 1$ to minimize the computation in each iteration. This test demonstrates the promise and capability of the SG-MDI algorithm for efficiently computing high-dimensional integrals.

| Dimension $(d)$ | Approximate Total nodes | $I_d(f)$ | | $I_d(\widehat{f})$ | |
|---|---|---|---|---|---|
| | | Relative error | CPU time(s) | Relative error | CPU time(s) |
| 10 | $1.0411 \times 10^6$ | $1.4725 \times 10^{-6}$ | 0.1541 | $2.0938 \times 10^{-5}$ | 0.1945 |
| 100 | $1.4971 \times 10^{60}$ | $1.5125 \times 10^{-5}$ | 80.1522 | $2.1235 \times 10^{-4}$ | 8.1920 |
| 300 | $3.3561 \times 10^{180}$ | $4.5377 \times 10^{-5}$ | 348.6000 | $6.3714 \times 10^{-4}$ | 52.0221 |
| 500 | $7.5230 \times 10^{300}$ | $7.7786 \times 10^{-5}$ | 1257.3354 | $1.0621 \times 10^{-3}$ | 219.8689 |
| 700 | $1.6767 \times 10^{421}$ | $1.0890 \times 10^{-4}$ | 3827.5210 | $1.4869 \times 10^{-3}$ | 574.9161 |
| 900 | $3.7524 \times 10^{541}$ | $1.4001 \times 10^{-4}$ | 9209.119 | $1.9117 \times 10^{-3}$ | 1201.65 |
| 1000 | $5.6136 \times 10^{601}$ | $1.5557 \times 10^{-4}$ | 13225.14 | $2.3248 \times 10^{-3}$ | 1660.84 |

Table 4.9: Computed results for $I_d(f)$ and $I_d(\widehat{f})$ by the SG-MDI algorithm.

## 5. Influence of parameters.

There are four parameters in the $d$-dimensional SG-MDI algorithm, they are respectively $r, m, s$, and $q$. Where $r \in \{1, 2, 3, 4\}$ represents the choice of one-dimensional numerical quadrature rule, namely, the (composite) trapezoidal rule ($r = 1$), Clenshaw-Curtis rule ($r = 2$), Gauss-Patterson rule ($r = 3$), and Gauss-Legendre rule ($r = 4$). The parameter $m$ stands for the multi-dimensional integration step size in the first stage of the algorithm so the dimension of the integration domain is reduced by $m$ in each iteration. In practice, $1 \leq m \leq q$, and it is preferred to be close to $q$. In this section we shall evaluate the performance of the SG-MDI algorithm when $m < q, m = q, m > q$. It should be noted that after $k := \lceil \frac{d}{m} \rceil$ iterations, the algorithm enters to its second stage and the iteration step size is changed to $s$. Since after the first stage, the remaining dimension $d - km$ small, so $s \in \{1, 2, 3\}$. It should be noted that after $k_1 := \lceil \frac{d - km}{s} \rceil$ iterations, the residual dimension satisfies $d - km - k_1 s \leq s$. Then in case $s = 2$ or $3$, one has two options to complete the algorithm. One either just continues the dimension reduction by calling 3d-SG-MDI or 2d-SG-MDI as explained in the definition of **Algorithm 2.3**, or compute the remaining 2- or 3-dimensional integral directly. The effect of these two choices will be tested in this section. Finally, the parameter $q$ represents the precision level of the sparse grid method. Obviously, the larger $q$ is, the higher accuracy of the computed results, the trade-off is more integration points must be used, hence, adds more cost. In this section we shall also test the impact of $q$ value to the SG-MDI algorithm.

### 5.1. Influence of parameter $r$.

We first examine the impact of one-dimensional quadrature rules, which are indicated by $r = 1, 2, 3, 4$, in the SG-MDI algorithm.

**Test 6:** Let $\Omega = [-1, 1]^d$ and choose the integrand $f$ as

$$(5.1) \qquad f(x) = \exp\Big(\sum_{i=1}^{d} (-1)^{i+1} x_i\Big), \qquad \widehat{f}(x) = \prod_{i=0}^{d} \frac{1}{0.9^2 + (x_i - 0.6)^2}.$$

Below we compare the performance of the SG-MDI algorithm with different $r$ in computing $I_d(f)$ and $I_d(\tilde{f})$ with the accuracy level $q = 10$ and step size $m = 10$.

| Parameter $(r)$ | Dimension $(d)$ | Approximate total nodes | Relative error | CPU time(s) |
|---|---|---|---|---|
| $r = 1$ | 10 | $1.5434 \times 10^7$ | $2.2076 \times 10^{-2}$ | 1.3309659 |
| | 30 | $3.6768 \times 10^{21}$ | $6.6228 \times 10^{-2}$ | 4.6347721 |
| | 50 | $8.7592 \times 10^{35}$ | $1.1038 \times 10^{-1}$ | 20.5783586 |
| | 70 | $2.0866 \times 10^{50}$ | $1.5453 \times 10^{-1}$ | 64.3242141 |
| | 90 | $4.9709 \times 10^{64}$ | $1.9868 \times 10^{-1}$ | 134.0323493 |
| | 100 | $7.6724 \times 10^{71}$ | $2.1035 \times 10^{-1}$ | 187.6477055 |
| $r = 2$ | 10 | $1.5434 \times 10^7$ | $3.2107 \times 10^{-5}$ | 0.5857483 |
| | 30 | $3.6768 \times 10^{21}$ | $9.6321 \times 10^{-5}$ | 3.0678083 |
| | 50 | $8.7592 \times 10^{35}$ | $1.6053 \times 10^{-4}$ | 12.2070886 |
| | 70 | $2.0866 \times 10^{50}$ | $2.2475 \times 10^{-4}$ | 34.9914719 |
| | 90 | $4.9709 \times 10^{64}$ | $2.8896 \times 10^{-4}$ | 81.3036934 |
| | 100 | $7.6724 \times 10^{71}$ | $3.3017 \times 10^{-4}$ | 113.0778901 |
| $r = 3$ | 10 | $1.0411 \times 10^6$ | $1.4725 \times 10^{-6}$ | 0.1540762 |
| | 30 | $1.1287 \times 10^{18}$ | $4.4177 \times 10^{-6}$ | 3.3439351 |
| | 50 | $1.2236 \times 10^{30}$ | $7.3628 \times 10^{-6}$ | 16.2212157 |
| | 70 | $1.3264 \times 10^{42}$ | $1.0308 \times 10^{-5}$ | 31.6057863 |
| | 90 | $1.4379 \times 10^{54}$ | $1.3253 \times 10^{-5}$ | 67.9195903 |
| | 100 | $1.4971 \times 10^{60}$ | $1.5125 \times 10^{-5}$ | 80.1522545 |
| $r = 4$ | 10 | $5.7789 \times 10^6$ | $2.2885 \times 10^{-6}$ | 0.1728588 |
| | 30 | $1.9299 \times 10^{20}$ | $6.8657 \times 10^{-6}$ | 6.2509345 |
| | 50 | $6.4454 \times 10^{33}$ | $1.1442 \times 10^{-5}$ | 31.420067 |
| | 70 | $2.1525 \times 10^{47}$ | $1.6020 \times 10^{-5}$ | 100.4530093 |
| | 90 | $7.1887 \times 10^{60}$ | $2.0597 \times 10^{-5}$ | 256.1780197 |
| | 100 | $4.1543 \times 10^{67}$ | $2.3115 \times 10^{-5}$ | 364.3777323 |

Table 5.1: Performance comparison of the SG-MDI algorithm with $q = 10, m = 10, s = 1$, and $r = 1, 2, 3, 4$.

Table 5.1 shows the computed results of $I(f)$ by the SG-MDI algorithm. We observe that the choice of one-dimensional quadrature rules has a significant impact on the accuracy and efficiency of the SG-MDI algorithm. The trapezoidal rule ($r = 1$) has the lowest precision and uses the most integration points, the Clenshaw-Curtis rule ($r = 2$) is the second lowest, and the Gauss-Patterson ($r = 3$) and Gauss-Legendre rule ($r = 4$) have the highest precision. Both the Clenshaw-Curtis and Gauss-Patterson rule use the nested grids, that is, the integration points of the ($q + 1$)th level contain those of the $q$th level. Although they use the same number of integration points, the Gauss-Patterson rule is more efficient than the Clenshaw-Curtis rule. Moreover, the Gauss-Patterson rule is more efficient than the Gauss-Legendre rule ($r = 4$) which uses the most CPU time and produces the most accurate solution. *This comparison suggests that the Gauss-Patterson rule is a winner among these four rules when they are used as the building blocks in the SG-MDI algorithm for high-dimensional integration.*

| Parameter $(r)$ | Dimension $(d)$ | Approximate total nodes | Relative error | CPU time(s) |
|---|---|---|---|---|
| $r = 1$ | 10 | $1.5434 \times 10^{7}$ | $3.6889 \times 10^{-2}$ | 0.9198790 |
| | 30 | $3.6768 \times 10^{21}$ | $1.1066 \times 10^{-1}$ | 2.7201514 |
| | 50 | $8.7592 \times 10^{35}$ | $1.8444 \times 10^{-1}$ | 4.6677602 |
| | 70 | $2.0866 \times 10^{50}$ | $2.5822 \times 10^{-1}$ | 6.7535072 |
| | 90 | $4.9709 \times 10^{64}$ | $3.3200 \times 10^{-1}$ | 9.2870665 |
| | 100 | $7.6724 \times 10^{71}$ | $3.5810 \times 10^{-1}$ | 10.5968167 |
| $r = 2$ | 10 | $1.5434 \times 10^{7}$ | $1.7211 \times 10^{-4}$ | 0.5905313 |
| | 30 | $3.6768 \times 10^{21}$ | $5.1633 \times 10^{-4}$ | 2.3826467 |
| | 50 | $8.7592 \times 10^{35}$ | $8.6056 \times 10^{-4}$ | 4.1569059 |
| | 70 | $2.0866 \times 10^{50}$ | $1.2047 \times 10^{-3}$ | 6.1272545 |
| | 90 | $4.9709 \times 10^{64}$ | $1.5490 \times 10^{-3}$ | 8.5287332 |
| | 100 | $7.6724 \times 10^{71}$ | $1.8013 \times 10^{-3}$ | 9.3035040 |
| $r = 3$ | 10 | $1.0411 \times 10^{6}$ | $2.0938 \times 10^{-5}$ | 0.1945813 |
| | 30 | $1.1287 \times 10^{18}$ | $6.2814 \times 10^{-5}$ | 0.7692118 |
| | 50 | $1.2236 \times 10^{30}$ | $1.0469 \times 10^{-4}$ | 1.8031252 |
| | 70 | $1.3264 \times 10^{42}$ | $1.4656 \times 10^{-4}$ | 3.6158604 |
| | 90 | $1.4379 \times 10^{54}$ | $1.8844 \times 10^{-4}$ | 6.2897488 |
| | 100 | $1.4971 \times 10^{60}$ | $2.1235 \times 10^{-4}$ | 8.1920274 |
| $r = 4$ | 10 | $5.7789 \times 10^{6}$ | $1.4304 \times 10^{-4}$ | 0.4950241 |
| | 30 | $1.9299 \times 10^{20}$ | $4.2912 \times 10^{-4}$ | 1.8407563 |
| | 50 | $6.4454 \times 10^{33}$ | $7.1521 \times 10^{-4}$ | 4.0708643 |
| | 70 | $2.1525 \times 10^{47}$ | $1.0012 \times 10^{-3}$ | 7.6650947 |
| | 90 | $7.1887 \times 10^{60}$ | $1.2873 \times 10^{-3}$ | 12.7835307 |
| | 100 | $4.1543 \times 10^{67}$ | $1.4410 \times 10^{-3}$ | 16.3940857 |

Table 5.2: Performance comparison of the SG-MDI algorithm with $q = 10, m = 10, s = 1$, and $r = 1, 2, 3, 4$.

Table 5.2 shows the computational results of $I_d(\widehat{f})$ by the SG-MDI algorithm. Similarly, the choice of one-dimensional quadrature rules has a significant impact on the accuracy and efficiency of the SG-MDI algorithm. Because the integrand $\widehat{f}(x)$ is simple, the SG-MDI algorithm with all four 1-d quadrature rules computes this integral very fast. Again, the trapezoidal rule is least accurate and the other three rules all perform very well, but a closer look shows that the Gauss-Patterson rule is again the best performer.

**5.2. Influence of parameter $m$.** From the Tables 4.5, 4.6, 5.1 and 5.2, we observe that when $m = 1$ is fixed, as the dimension $d$ increases, the number of iterations by the SG-MDI algorithm also increases, and the computational efficiency decreases rapidly. In practice, the step size $m$ of the SG-MDI algorithm in the first stage iteration should not be too large or too small. One strategy is to use variable step sizes. After selecting an appropriate initial step size $m$, it can be decreased during the dimension iteration. The next test presents a performance comparison of the SG-MDI algorithm for $m < q, m = q, m > q$.

**Test 7.** Let $\Omega = [-1, 1]^d$, $f$ and $\widehat{f}$ be the same as in (4.4).

We compute these integrals using the SG-MDI algorithm with $s = 1, r = 3$ (Gauss-Patterson rule) and $q = 10$.

| Parameter $(r)$ | Dimension $(d)$ | number of iterations | Relative error | CPU time(s) |
|---|---|---|---|---|
| $m = 5$ | 10 | 2 | $2.7427 \times 10^{-6}$ | 0.6076211 |
| | 30 | 6 | $8.2283 \times 10^{-6}$ | 2.4312369 |
| | 50 | 10 | $1.3713 \times 10^{-5}$ | 5.8717692 |
| | 70 | 14 | $1.9199 \times 10^{-5}$ | 11.9989244 |
| | 90 | 18 | $2.4684 \times 10^{-5}$ | 21.2031731 |
| | 100 | 20 | $2.8225 \times 10^{-5}$ | 27.3314984 |
| $m = 10$ | 10 | 1 | $2.0938 \times 10^{-5}$ | 0.1945813 |
| | 30 | 3 | $6.2814 \times 10^{-5}$ | 0.7692118 |
| | 50 | 5 | $1.0469 \times 10^{-4}$ | 1.8031252 |
| | 70 | 7 | $1.4656 \times 10^{-4}$ | 3.6158604 |
| | 90 | 9 | $1.8844 \times 10^{-4}$ | 6.2897488 |
| | 100 | 10 | $2.1235 \times 10^{-4}$ | 8.1920274 |
| $m = 15$ | 10 | 1 | $2.0938 \times 10^{-5}$ | 0.1945813 |
| | 30 | 2 | $1.1146 \times 10^{-3}$ | 1.0704125 |
| | 50 | 4 | $1.8577 \times 10^{-3}$ | 1.9306281 |
| | 70 | 5 | $2.6007 \times 10^{-3}$ | 3.1370912 |
| | 90 | 6 | $3.3438 \times 10^{-3}$ | 4.5365430 |
| | 100 | 7 | $3.7154 \times 10^{-3}$ | 5.6229055 |

Table 5.3: Efficiency comparison of the SG-MDI algorithm with $q = 10$, $r = 3$, $s = 1$ and $m = 5, 10, 15$.

Table 5.3 presents the computed results of integral $I(f)$ in **Test 7** by the SG-MDI algorithm. It is easy to see that the accuracy and efficiency of the SG-MDI algorithm with different $m$ are different, this is because the step size $m$ affects the number of iterations in the first stage. It shows that the larger step size $m$, the smaller the number of iterations and the number of symbolic functions that need to be saved, and the fewer CPU time are used, but at the expense of decreasing accuracy and higher truncation error. On the other hand, the smaller step size $m$, the more accurate of the computed results, but at the expense of more CPU time. To explain this observation, we notice that when the step size $m$ is large, each iteration reduces more dimensions, and the number of symbolic functions that need to be saved is less, but the truncation error generated by the system is larger, although the used CPU time is small. When the step size is small, however, more symbolic functions need to be saved. This is because each iteration reduces a small number of dimensions, the error of computed result is smaller, but the used CPU time is larger. We also observed that the SG-MDI algorithm with step size $m \approx p$ achieves a good balance between CPU time and accuracy.

Table 5.4 shows the computed results of $I_d(\widehat{f})$ in **Test 7** the SG-MDI algorithm. As expected, choosing different parameters $m$ has a significant impact on the accuracy and efficiency of the SG-MDI algorithm. Since function evaluation of the integrand $\widehat{f}(x)$ is more complicated, the influence of different parameters $m$ on the SG-MDI algorithm is dramatic. Again, the SG-MDI algorithm with step size $m \approx p$ strikes a balance between CPU time and accuracy.

| Parameter $(r)$ | Dimension $(d)$ | number of iterations | Relative error | CPU time(s) |
|---|---|---|---|---|
| $m = 5$ | 10 | 2 | $3.5104 \times 10^{-7}$ | 0.7389749 |
| | 30 | 6 | $1.0531 \times 10^{-6}$ | 45.8154739 |
| | 50 | 10 | $1.7552 \times 10^{-6}$ | 418.3180822 |
| | 70 | 14 | $2.4573 \times 10^{-6}$ | 2063.0159854 |
| | 90 | 18 | $3.1594 \times 10^{-6}$ | 6700.8943929 |
| | 100 | 20 | $3.6207 \times 10^{-6}$ | 12916.407730 |
| $m = 10$ | 10 | 1 | $1.6855 \times 10^{-6}$ | 0.295855 |
| | 30 | 3 | $5.0567 \times 10^{-6}$ | 29.5098187 |
| | 50 | 5 | $8.4279 \times 10^{-6}$ | 293.1338336 |
| | 70 | 7 | $1.1799 \times 10^{-5}$ | 1464.5795716 |
| | 90 | 9 | $1.5170 \times 10^{-5}$ | 2388.3820738 |
| | 100 | 10 | $1.7825 \times 10^{-5}$ | 3077.1890005 |
| $m = 15$ | 10 | 1 | $1.6855 \times 10^{-6}$ | 0.295855 |
| | 30 | 2 | $2.4796 \times 10^{-4}$ | 41.6601228 |
| | 50 | 4 | $4.1327 \times 10^{-4}$ | 159.7405442 |
| | 70 | 5 | $5.7859 \times 10^{-4}$ | 432.4688346 |
| | 90 | 6 | $7.4390 \times 10^{-4}$ | 861.6527617 |
| | 100 | 7 | $8.2655 \times 10^{-4}$ | 1153.8207432 |

Table 5.4: Efficiency comparison of the SG-MDI algorithm with $q = 10$, $r = 3$, $s = 1$ and $m = 5, 10, 15$.

**5.3. Influence of the parameter** $s$**.** In this subsection, we test the impact of the step size $s$ used in the second stage of the SG-MDI algorithm. Recall that the step size $m \approx p$ works well in the first stage. After $k := \lceil \frac{d}{m} \rceil$ iterations, the dimension of the integration domain is reduced to $d - km$ which is relatively small, hence, $s$ should be small. In practice, we have $1 \le s \le 3$. The goal of the next test is to provide a performance comparison of the SG-MDI algorithm with $s = 1, 2, 3$.

**Test 8.** Let $\Omega = [-1, 1]^d$, and choose the integrand $f$ as

$$(5.2) \qquad f(x) = \prod_{i=0}^{d} \frac{1}{0.9^2 + (x_i - 0.6)^2}.$$

We compute this integral $I(f)$ using the SG-MDI algorithm with $m = 10, r = 3$ (Gauss-Patterson rule) and $p = 10$. Table 5.5 displays the computed results. We observe that the same accuracy is achieved in all cases $s = 1, 2, 3$, which is expected. Moreover, the choice of $s$ has little effect on the efficiency of the algorithm. An explanation for this observation is that because $d - km$ becomes small after the first stage iterations, so the number of the second stage iterations is small, the special way of performing function evaluations in the SG-MDI algorithm is not sensitive to the small variations in the choice of $s = 1, 2, 3$.

| Dimension $(d)$ | Approximate total nodes | Relative error | $s = 1$ CPU time(s) | $s = 2$ CPU time(s) | $s = 3$ CPU time(s) |
|---|---|---|---|---|---|
| 10 | $1.0411 \times 10^6$ | $2.0938 \times 10^{-5}$ | 0.19458 | 0.212757 | 0.249253 |
| 30 | $1.1287 \times 10^{18}$ | $6.2814 \times 10^{-5}$ | 0.76921 | 1.013215 | 1.056075 |
| 50 | $1.2236 \times 10^{30}$ | $1.0469 \times 10^{-4}$ | 1.80312 | 1.964943 | 2.088460 |
| 70 | $1.3264 \times 10^{42}$ | $1.4656 \times 10^{-4}$ | 3.61586 | 4.080034 | 3.931925 |
| 90 | $1.4379 \times 10^{54}$ | $1.8844 \times 10^{-4}$ | 6.28975 | 6.809549 | 6.916479 |
| 100 | $1.4971 \times 10^{60}$ | $2.1235 \times 10^{-4}$ | 8.19202 | 8.621845 | 8.587019 |

Table 5.5: Efficiency comparison of the SG-MDI algorithm with $r = 3$ and $s = 1, 2, 3$.

**5.4. Influence of the parameter $p$ or $N$.** Finally, we examine the impact of the accuracy level $p$ on the SG-MDI algorithm. Recall that the parameter $p$ in the SG method is related to the number of integration points $N$ in one coordinate direction on the boundary. It is easy to check that for the trapezoidal ($r = 1$) and Clenshaw Curtis ($r = 2$) quadrature rules, $q = 1, 2, 3, 4, 5, 6, 7, 8, 9, 10$ correspond to $N = 1, 3, 5, 9, 9, 17, 17, 17, 17, 33$, for the Gauss-Patterson quadrature rule ($r = 3$), $q = 1, 2, 3, 4, 5, 6, 7, 8, 9, 10$ correspond to $N = 1, 3, 3, 7, 7, 7, 15, 15, 15, 15$, and for the Gauss-Legendre quadrature rule ($r = 4$), $q = 1, 2, 3, 4, 5, 6, 7, 8, 9, 10, 12, 14$ correspond to $N = 1, 2, 3, 4, 5, 6, 7, 8, 9, 10, 12, 14$. Therefore, we only need to examine the impact of the parameter $N$ on the SG-MDI algorithm. To the end, we consider the case $m = 1, s = 1$, and $r = 4$ (Gauss-Legendre rule) in the next test.

**Test 8.** Let $\Omega = [-1, 1]^d$ and choose the following integrands:

$$f(x) = \exp\Big(\sum_{i=1}^{d}(-1)^{i+1}x_i\Big), \quad \widehat{f}(x) = \cos\Big(2\pi + \sum_{i=1}^{d}x_i\Big), \quad \widetilde{f}(x) = \prod_{i=0}^{d}\frac{1}{0.9^2 + (x_i - 0.6)^2}.$$

Table 5.6, 5.7, and 5.8 present respectively the computed results for $d = 5, 10$ and $N = 4, 6, 8, 10, 12, 14$. We observe that the quality of approximations also depends on the behavior of the integrand. For very oscillatory and rapidly growing functions, more integration points must be used to achieve good accuracy.

**6. Computational complexity.**

**6.1. The relationship between the CPU time and $N$.** In this subsection, we examine the relationship between the CPU time and parameter $N$ using the regression technique based on the test data.

**Test 9.** Let $\Omega, f, \widehat{f}$, and $\widetilde{f}$ be the same as in **Test 8**.

Figure 6.1 and 6.2 show the CPU time as a function of $N$ obtained by the least square method and the fitting functions are given in Table 6.1. All the results indicate that the CPU time grows at most linearly in $N$.

| $p(N)$ | $d = 5$ | | | $d = 10$ | | |
|---|---|---|---|---|---|---|
| | Approximate total nodes | Relative error | CPU time(s) | Approximate total nodes | Relative error | CPU time(s) |
| 4(4) | 241 | $3.8402 \times 10^{-3}$ | 0.1260 | 1581 | $5.7516 \times 10^{-2}$ | 0.3185 |
| 6(6) | 2203 | $1.6039 \times 10^{-5}$ | 0.1608 | 40405 | $2.3524 \times 10^{-3}$ | 0.4546 |
| 8(8) | 13073 | $1.7195 \times 10^{-8}$ | 0.2127 | 581385 | $3.6774 \times 10^{-5}$ | 0.6056 |
| 10(10) | 58923 | $6.4718 \times 10^{-12}$ | 0.2753 | 5778965 | $2.2885 \times 10^{-7}$ | 0.7479 |
| 12(12) | 218193 | $1.8475 \times 10^{-12}$ | 0.3402 | 44097173 | $2.2746 \times 10^{-9}$ | 1.0236 |
| 14(14) | 695083 | $8.0013 \times 10^{-12}$ | 0.4421 | 112613833 | $3.8894 \times 10^{-11}$ | 1.2377 |

Table 5.6: Performance comparison of the SD-MDI algorithm with $p, N = 4, 6, 8, 10, 12, 14$ for computing $I_d(f)$.

| $p(N)$ | $d = 5$ | | | $d = 10$ | | |
|---|---|---|---|---|---|---|
| | Approximate total nodes | Relative error | CPU time(s) | Approximate total nodes | Relative error | CPU time(s) |
| 4(4) | 241 | $1.4290 \times 10^{-2}$ | 0.1664 | 1581 | $8.1129 \times 10^{-1}$ | 0.3174 |
| 6(6) | 2203 | $6.1319 \times 10^{-5}$ | 0.2159 | 40405 | $3.6823 \times 10^{-2}$ | 0.4457 |
| 8(8) | 13073 | $6.6347 \times 10^{-8}$ | 0.2526 | 581385 | $5.8931 \times 10^{-4}$ | 0.5571 |
| 10(10) | 58923 | $2.5247 \times 10^{-11}$ | 0.3305 | 5778965 | $3.8749 \times 10^{-6}$ | 0.6717 |
| 12(12) | 218193 | $1.7163 \times 10^{-12}$ | 0.3965 | 44097173 | $2.2490 \times 10^{-8}$ | 0.8843 |
| 14(14) | 695083 | $8.2889 \times 10^{-12}$ | 0.5277 | 112613833 | $8.7992 \times 10^{-10}$ | 1.1182 |

Table 5.7: Performance comparison of the SD-MDI algorithm with $p, N = 4, 6, 8, 10, 12, 14$ for computing $I_d(\hat{f})$.

| $p(N)$ | $d = 5$ | | | $d = 10$ | | |
|---|---|---|---|---|---|---|
| | Approximate total nodes | Relative error | CPU time(s) | Approximate total nodes | Relative error | CPU time(s) |
| 4(4) | 241 | $6.1894 \times 10^{-4}$ | 0.1275 | 1581 | $1.5564 \times 10^{-2}$ | 0.1657 |
| 6(6) | 2203 | $1.9354 \times 10^{-3}$ | 0.1579 | 40405 | $1.0163 \times 10^{-2}$ | 0.2530 |
| 8(8) | 13073 | $1.5488 \times 10^{-4}$ | 0.1755 | 581385 | $2.2076 \times 10^{-3}$ | 0.3086 |
| 10(10) | 58923 | $1.7878 \times 10^{-6}$ | 0.1963 | 5778965 | $1.4304 \times 10^{-4}$ | 0.3889 |
| 12(12) | 218193 | $7.0609 \times 10^{-7}$ | 0.2189 | 44097173 | $9.3339 \times 10^{-6}$ | 0.4493 |
| 14(14) | 695083 | $1.7194 \times 10^{-8}$ | 0.2459 | 112613833 | $2.4671 \times 10^{-7}$ | 0.4864 |

Table 5.8: Performance comparison of the SD-MDI algorithm with $p, N = 4, 6, 8, 10, 12, 14$ for computing $I_d(\tilde{f})$.

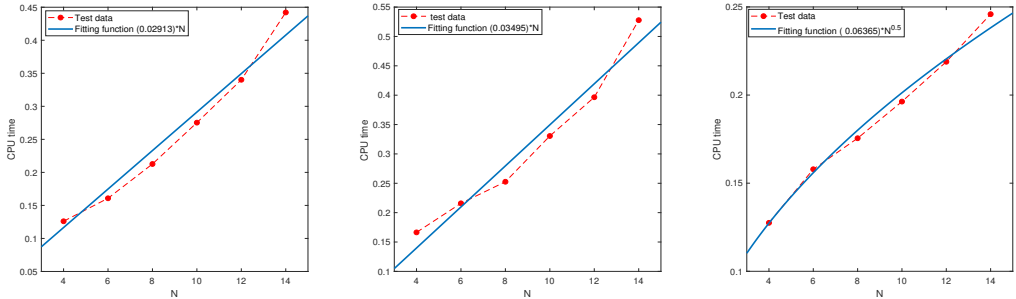

Fig. 6.1: The relationship between the CPU time and parameter $N$ when $d = 5$ for computing $I_d(f)$ (left), $I_d(\hat{f})$ (middle), $I_d(\tilde{f})$ (right).

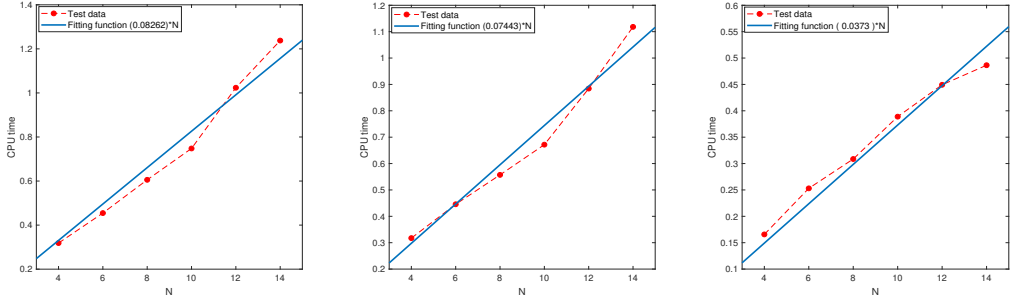

Fig. 6.2: The relationship between the CPU time and parameter $N$ when $d = 10$ for computing $I_d(f)$ (left), $I_d(\widehat{f})$ (middle), $I_d(\widetilde{f})$ (right).

| Integrand | $r$ | $m$ | $d$ | Fitting function | $R$-square |
|-----------|-----|-----|-----|------------------|------------|
| $f(x)$ | 3 | 1 | 5 | $h_1(N) = (0.02913) * N$ | 0.9683 |
| $\widehat{f}(x)$ | 3 | 1 | 5 | $h_2(N) = (0.03495) * N$ | 0.9564 |
| $\widetilde{f}(x)$ | 3 | 1 | 5 | $h_3(N) = (0.06365) * N^{\frac{1}{2}}$ | 0.9877 |
| $f(x)$ | 3 | 1 | 10 | $h_4(N) = (0.08262) * N$ | 0.9692 |
| $\widehat{f}(x)$ | 3 | 1 | 10 | $h_5(N) = (0.07443) * N$ | 0.9700 |
| $\widetilde{f}(x)$ | 3 | 1 | 10 | $h_6(N) = (0.0373) * N$ | 0.9630 |

Table 6.1: The relationship between the CPU time and parameter $N$.

## 6.2. The relationship between the CPU time and the dimension $d$.
Recall that the computational complexity of the sparse grid method is of the order $O(N \cdot (\log N)^{d-1})$ for computing $I_d(f)$, which grows exponentially in $d$ with base $\log N$. The numerical tests presented above overwhelmingly and consistently indicate that the SG-MDI algorithm has hidden capability to overcome the curse of dimensionality which hampers the sparse grid method. The goal of this subsection is to find out the computational complexity of the SG-MDI algorithm (in terms of CPU time as a function of $d$) using the regression technique based on numerical test data.

**Test 10.** Let $\Omega = [0, 1]^d$ and consider the following five integrands:

$$f_1(x) = \exp\Big(\sum_{i=1}^{d}(-1)^{i+1}x_i\Big), \qquad f_2(x) = \prod_{i=0}^{d}\frac{1}{0.9^2 + (x_i - 0.6)^2},$$

$$f_3(x) = \frac{1}{\sqrt{2\pi}}\exp\Big(-\frac{1}{2}|x|^2\Big), \qquad f_4(x) = \cos\Big(2\pi + \sum_{i=1}^{d}x_i\Big),$$

$$f_5(x) = \frac{1}{2^d}\exp\Big(\sum_{i=1}^{d}(-1)^{i+1}x_i\Big).$$

Figure 6.3 displays the CPU time as functions of $d$ obtained by the least square regression method whose analytical expressions are given in Table 6.2. We note that the parameters of the SG-MDI algorithm only affect the coefficients of the fitting functions, but not their orders in $d$.

| Integrand | $r$ | $m$ | $s$ | $q(N)$ | Fitting function | $R$-square |
|-----------|-----|-----|-----|--------|------------------|------------|
| $f_1$ | 1 | 10 | 1 | 10(33) | $g_1 = (1.092e-05) * Nd^3$ | 0.9995 |
|       | 2 | 10 | 1 | 10(33) | $g_2 = (6.531e-06) * Nd^3$ | 0.9977 |
|       | 3 | 10 | 1 | 10(15) | $g_3 = (8.076e-05) * Nd^{2.4}$ | 0.9946 |
|       | 4 | 10 | 1 | 10(10) | $g_4 = (3.461e-05) * Nd^3$ | 0.9892 |
| $f_2$ | 1 | 10 | 1 | 10(33) | $g_5 = 0.003820 * Nd^{1.1}$ | 0.9985 |
|       | 2 | 10 | 1 | 10(33) | $g_6 = 0.003432 * Nd^{1.1}$ | 0.9986 |
|       | 3 | 5 | 1 | 10(15) | $g_6 = (7.152e-05) * Nd^{2.2}$ | 0.9983 |
|       | 3 | 10 | 1 | 10(15) | $g_9 = (1.106e-07) * Nd^3$ | 0.9998 |
|       | 3 | 15 | 1 | 10(15) | $g_7 = 0.0004145 * Nd^{1.47}$ | 0.9955 |
|       | 3 | 15 | 2 | 10(15) | $g_7 = (5.681e-05) * Nd^2$ | 0.9961 |
|       | 3 | 15 | 3 | 10(15) | $g_7 = (5.677e-05) * Nd^2$ | 0.9962 |
|       | 4 | 10 | 1 | 10(10) | $g_8 = 0.00016 * Nd^2$ | 0.9965 |
| $f_3$ | 3 | 5 | 1 | 10(15) | $g_{11} = (8.312e-08) * Nd^3$ | 0.9977 |
|       | 3 | 10 | 1 | 10(15) | $g_{12} = 0.0008441 * Nd^{2.7}$ | 0.9844 |
|       | 3 | 15 | 1 | 10(15) | $g_{13} = 0.0003023 * Nd^{2.8}$ | 0.9997 |
| $f_4$ | 3 | 10 | 1 | 10(15) | $g_{18} = (4.053e-05) * Nd^3$ | 0.9903 |
| $f_5$ | 3 | 10 | 1 | 10(15) | $g_{19} = (8.461e-07) * Nd^3$ | 0.9958 |

Table 6.2: The relationship between CPU time as a function of the dimension $d$.

We also quantitatively characterize the performance of the fitted curve by the $R$-square in Matlab, which is defined as $R\text{-square} = 1 - \frac{\sum_i^n (y_i - \widehat{y}_i)^2}{\sum_i^n (y_i - \overline{y})^2}$. Where $y_i$ represents a test data output, $\widehat{y}_i$ refers to the predicted value, and $\overline{y}$ indicates the mean value of $y_i$. Table 6.2 also shows that the $R$-square of all fitting functions is very close to 1, which indicates that the fitting functions are quite accurate. These results suggest that the CPU time grows at most cubically in $d$. Combining the results of **Test 8** in Section 5.4 we conclude that the CPU time required by the proposed SD-MDI algorithm grows at most in the polynomial order $O(d^3 N)$.

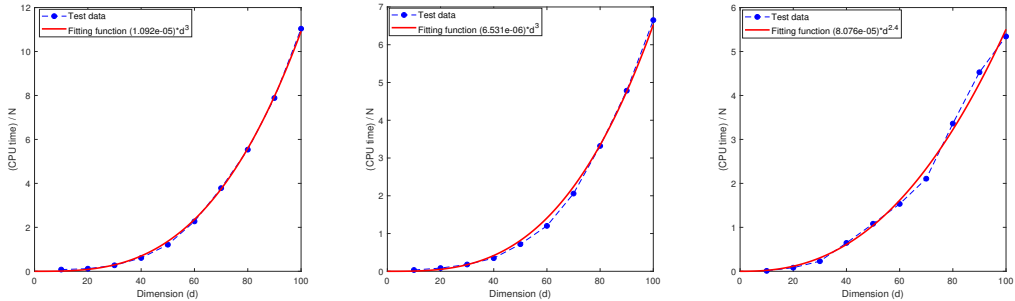

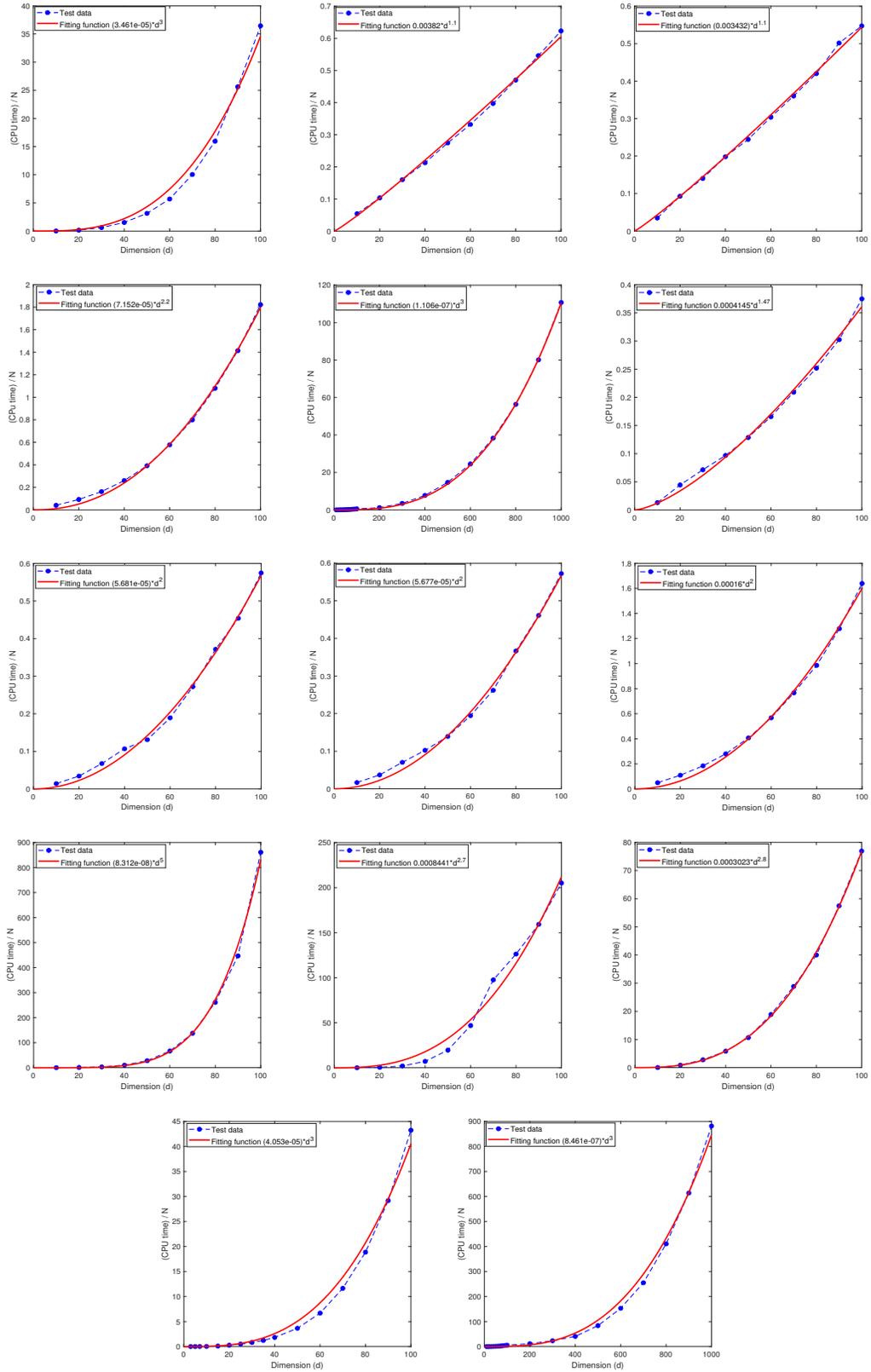

Fig. 6.3: The relationship between the CPU time and dimension $d$.

**7. Conclusions.** This paper presented an efficient and fast implementation algorithm (or solver), called the SD-MDI algorithm, for high-dimensional numerical integration using the sparse grid method. It is based on combining the idea of dimension iteration/reduction combined with the idea of computing the function evaluations at all integration points in cluster so many computations can be reused. It was showed numerically that the computational complexity (in terms of the CPU time) of the SD-MDI algorithm grows at most cubically in the dimension $d$, and overall in the order $O(Nd^3)$, where $N$ denotes the maximum number of integration points in each coordinate direction. This shows that the SD-MDI algorithm could effectively circumvent the curse of the dimensionality in high-dimensional numerical integration, hence, makes sparse grid methods not only become competitive but also can excel for the job. Extensive numerical tests were provided to examine the performance of the SD-MDI algorithm and to carry out performance comparisons with the standard sparse grid method and with the Monte Carlo (MC) method. It was showed that the SD-MDI algorithm (regardless the choice of the 1-d base sparse grid quadrature rules) is faster than the MC method in low and medium dimensions (i.e., $d \lesssim 100$), much faster in very high dimensions (i.e., $d \approx 1000$), and succeeds even when the MC method fails. An immediate application of the proposed SD-MDI algorithm is to solve high-dimensional PDEs which will be reported in a forthcoming work.